\documentclass[a4paper, 11pt]{amsart}
\usepackage[a4paper, margin=2.5cm]{geometry}
\usepackage{amsmath}
\usepackage{amssymb}
\usepackage{amsfonts}
\usepackage{amscd}
\usepackage{amsthm}
\usepackage{mathrsfs}
\usepackage[bb=boondox,bbscaled=.95,cal=boondoxo]{mathalfa}
\usepackage[numbers]{natbib}

\newcommand{\CC}{{\rm\bf C}}
\newcommand{\RR}{{\rm\bf R}}

\newcommand{\ZZ}{{\rm\bf Z}}
\newcommand{\NN}{{\rm\bf N}}

\newcommand{\PP}{{\rm\bf P}}

\DeclareMathOperator{\Spec}{\mathrm{Spec}}
\DeclareMathOperator{\diag}{\mathrm{diag}}

\DeclareMathOperator{\GL}{\mathrm{GL}}
\DeclareMathOperator{\SL}{\mathrm{SL}}

\DeclareMathOperator{\SO}{\mathrm{SO}}
\DeclareMathOperator{\SU}{\mathrm{SU}}

\DeclareMathOperator{\Sym}{\mathrm{Sym}}

\DeclareMathOperator{\Hom}{\mathrm{Hom}}
\DeclareMathOperator{\Ext}{\mathrm{Ext}}

\DeclareMathOperator{\AS}{{\mathrm{AS}}}

\DeclareMathOperator{\Ind}{\mathrm{Ind}}

\DeclareMathOperator{\Lie}{\mathrm{Lie}}

\newcommand{\liea}{{\mathfrak {a}}}
\newcommand{\lieb}{{\mathfrak {b}}}
\newcommand{\liec}{{\mathfrak {c}}}

\newcommand{\lieg}{{\mathfrak {g}}}
\newcommand{\lieh}{{\mathfrak {h}}}
\newcommand{\liek}{{\mathfrak {k}}}

\newcommand{\lien}{{\mathfrak {n}}}

\newcommand{\liet}{{\mathfrak {t}}}

\newcommand{\liep}{{\mathfrak {p}}}

\newcommand{\liesl}{{\mathfrak {sl}}}

\theoremstyle{plain}
\newtheorem{theorem}{Theorem}[section]

\newtheorem{proposition}[theorem]{Proposition}

\theoremstyle{remark}
\newtheorem{remark}[theorem]{Remark}

\bibpunct{[}{]}{,}{n}{}{,}

\begin{document}

\title{Towards a tensor-classification of Harish-Chandra modules:\\The case $\SL_2(\RR)$}

\author{Fabian Januszewski}
\address{Institut f\"ur Mathematik, Fakult\"at EIM, Paderborn University, Germany}
\email{fabian.januszewski@math.upb.de}
        
\maketitle

\begin{abstract}
  We consider the category of Harish-Chandra modules for $\SL_2(\RR)$ as a module over the category of finite-dimensional representations of $\SL(2)$ with respect to the tensor product. In this note we use classical results about principal series to obtain a classification of $\otimes$-submodules of the category of Harish-Chandra modules of $\SL_2(\RR)$. The resulting classification recovers the classical classification of irreducible Harish-Chandra modules. Our methods are expected to generalize to arbitrary reductive pairs and more general base fields.
\end{abstract}


\section*{Introduction}

The study of (unitary) representations of reductive Lie groups has its origins in the pioneering work of Elie Cartan and Hermann Weyl yielding the classification of irreducible unitary representations of (connected) compact Lie groups. In this case, all unitary representations are semi-simple and irreducibles turn out to be finite-dimensional.

The non-compact case is significantly more challenging: infinite-dimensional irreducible representations enter the scene and the category of finite length modules is no more semi-simple. The irreducible representation of $\SL_2(\RR)$ were classified by Bargmann \cite{bargmann1947} and independently by Gelfand and Naimark \cite{gelfandnaimark1946}. After proving the Plancherel formula for $\SL_2(\RR)$ in \cite{harishchandra1953a}, Harish-Chandra subsequently laid the foundation for the classification of irreducible admissible representations for general reductive Lie groups in \cite{harishchandra1953b,harishchandra1954a,harishchandra1954b} which in turn led Langlands to formulate his classification \cite{book_langlands1989} that now bears his name.

\smallskip
Harish-Chandra replaced a reductive Lie group $G$ (in Harish-Chandra's class) by a reductive pair $({\mathfrak{g}},K)$, where ${\mathfrak{g}}$ is the complexified Lie algebra of $G$ and $K$ is a maximal compact Lie subgroup of $G$ and reduced the problem of classifying unitary representations of $G$ to classifying infinitesimally unitary representations of $({\mathfrak{g}},K)$. In Harish-Chandra's work both algebra and analysis have their rightful places.

\smallskip
The definition of a $({\mathfrak{g}},K)$-module is a purely algebraic concept and is meaningful over arbitrary fields of characteristic $0$, cf.~\cite{harder2014,harris2012,januszewski2018}. This point of view was initially motivated by arithmetic applications (cf.~\cite{harder2014,harris2012,harris2012erratum,januszewski2018,januszewski2019} for instance) and has been also adopted in the context of contractions of representations in \cite{bernsteinhigsonsubag2020}.

Therefore, we are naturally led at complementing and extending the existing algebraic and homological methods in the theory of Harish-Chandra modules in order to obtain a self-contained algebraic theory independent of the ground field (or even ground ring).

\smallskip
With Beilinson-Bernstein's classification of $({\mathfrak{g}},K)$-modules in terms of infinitesimal characters and sheaves on $K$-orbits on the flag variety of $G$ (provided $G$ originates from a linear algebraic group) we have a clean algebraic theory over algebraically closed fields of characteristic $0$ at hand, which currently is being extended to general bases in joint work of the author with Takuma Hayashi \cite{januszewskihayashi-preprint}.

\smallskip
In this note, we embark on a categorical approach to the classification problem, motivated by the following observations.

\smallskip
When restricting to finite-dimensional (rational) representations of an algebraic group, the theory of Tannakian categories implies that the $\otimes$-structure together with a suitable fibre functor is sufficient to recover the unterlying group. Therefore, the same structure is morally sufficient to classify finite-dimensional irreducible representations.

Although the concept of Tannakian categories does not apply to situation at hand which involves infinite-dimensional representations of a single reductive pair, with the notion of coherent families and translation functors the tensor product always had its rightful place in the theory of Harish-Chandra modules.

On a conceptual level, the $\otimes$-structure on $({\mathfrak{g}},K)$-modules plays a central role in the author's homological character theory \cite{januszewski2014}: In order to obtain a well defined and well behaved theory, a suitable Weyl denominator has to be inverted with respect to a multiplication induced by the tensor product. Now since this character theory classifies irreducible $({\mathfrak{g}},K)$-modules, this suggests that the $\otimes$-structure plays a central role in the classification problem.

\smallskip
In this note, we put the $\otimes$-structure at the center of attention.

\smallskip
Balmer \cite{balmer2005} attaches to any small $\otimes$-triangulated category $\mathcal D$ a topological space ${\mathrm{Spec}}\,\mathcal D$ which is universal for a suitable notion of support. In particular, this allows Balmer to define a notion of support for objects (and $\otimes$-ideals) in the derived category under consideration.

Aiming at a (complete) description of the Balmer spectrum of nice categories of $({\mathfrak{g}},K)$-modules, such as finite length modules, admissible modules, etc., therefore appears as a natural problem, which in other settings has been successfully achieved (Hopkins-Neeman's classification of localising subcategories of the derived category of modules over a commutative ring \cite{book_hopkins1987,neeman1992}, Benson-Iyengar-Krause classified localising subcagetories of the stable category of $k[G]$-modules for a finite group $G$ \cite{bensoniyengarkrause2008}, extending the Benson-Carlson-Rickard classification \cite{bensoncarlonsrickard1997}, and other classification results, cf.\ \cite{balmer2010}).

\smallskip
There is a technical obstacle however: Neither the category of finite length modules nor the category of admissible modules are closed under the tensor product. The observation that in order to obtain a meaningful algebraic / homological character theory as in \cite{januszewski2014}, the only relevant tensor products are the ones with finite-dimensional modules, we restrict our attention to the following setup:

\smallskip
We consider the category of finite length $({\mathfrak{g}},K)$-modules as a {\em module} over the category of finite-dimensional $({\mathfrak{g}},K)$-modules with respect to the tensor product (in the sense of Yetter \cite{yetter1992,yetter2003}) and aim at classifying the all thick $\otimes$-submodules of said category.

\smallskip
Stevenson \cite{stevenson2013} defined a notion of support in such situations on the Balmer spectrum of the smaller category operating via a suitable operation on a larger category. However, in our case the Balmer spectrum of the category of finite-dimensional modules is trivial, because the underlying category is semi-simple.

\smallskip
Nonetheless we demonstrate that already in the first non-trivial case $G=\SL_2(\RR)$, the classification of thick $\otimes$-submodules of the category of finite-length modules over the category of finite-dimensional modules is meaningful and reflects the classical classification of irreducible Harish-Chandra modules in this case. In fact, together with the infinitesimal character, it is possible to deduce a complete classification of irreducible modules (cf.\ section \ref{sec:completeclassification}). Furthermore, there is a close relationship between the classification of thick $\otimes$-submodules and $K$-orbits on the flag variety of $G$ (cf.\ section \ref{sec:geometricinterpretation}).

\smallskip
For an expert in the representation theory of real reductive groups neither the classification in Theorem \ref{thm:classification}, nor the techniques involved in its proof will be surprising. Nonetheless we included many details in order to make this text accessible to non-experts as well.

\subsubsection*{Outline of the paper}

This note is organized as follows. In section \ref{sec:generalities} we introduce notation and collect relevant technical results. Since all those statements are known in the literature for general reductive pairs $({\mathfrak{g}},K)$, there is no reason to restrict the statements in section \ref{sec:generalities} to a special case.

\smallskip
In section \ref{sec:categorical} we introduce the categorical classification problem for general reductive pairs $(\lieg,K)$.

\smallskip
From section \ref{sec:sl2} on we restrict ourselves to $G=\SL_2(\RR)$ and after a short discussion of finite-dimensional representations, we focus on composition series of principal series representations, which are one key ingredient to our classification: Composition factors of tensor products of principal series representations with finite-dimensional representations are easily decomposed in the case at hand. The other key ingredient towards the classification is provided by the asymptotic $K$-support, the latter being stable under tensor products.

\smallskip
At the end of this section \ref{sec:sl2} we summarize the classification in Theorem \ref{thm:classification} and discuss the relation to Beilinson-Bernstein localization and the classification problem for irreducible modules.

\subsubsection*{Acknowledgement}

The authors thanks Henning Krause for stimulating conversations and for references for known classifications in $\otimes$-triangulated contexts.

\section{Generalities}\label{sec:generalities}

\subsection{Reductive pair}

In this section $G$ denotes a connected linear reductive group over $\RR$. We write $G_\CC$ for its complexification and $G(\RR)$ for the group of real-valued points. Then $G(\RR)$ is a real reductive Lie group. Fix a maximal compact subgroup of $K$ of $G(\RR)$ and write $\lieg=\Lie(G)\otimes_\RR\CC$ for the complexified Lie algebra. Then $(\lieg,K)$ admits the structure of a {\em reductive pair} in the sense of Chap.\ IV in \cite{book_knappvogan1995}. We write $\theta:\lieg\to\lieg$ for the Cartan involution attached to $K$ and denote by $\lieg=\liep\oplus\liek$ the Cartan decomposition accordingly. As the notation suggests, $\liek$ denotes the Lie algebra of $K$ (or rather of the connected compontent of the identity in $K$).

Decompose $\lieg=\lien^-\oplus\lieh\oplus\lien$ for a Cartan subalgebra $\lieh\subseteq\lieg$ and a choice of positive root spaces $\lien$ and its opposite $\lien^-$ ($\lien$ corresponds to a choice of Borel subalgebra $\lieb=\lieh+\lien$ inside $\lieg$).

\subsection{Subpairs}

A subpair $(\liea,B)$ of $(\lieg,K)$ consists of a Lie-subalgebra $\liea$ of $\lieg$ and a closed subgroup $B\subseteq K$ such that its Lie algebra $\lieb\subseteq\liea$ is contained in $\liea$. Subpairs of subpairs are defined similarly and saying that the pair $(\liea,B)$ is $\theta$-stable is equivalent to the condition $\theta(\liea)=\liea$.

\subsection{Universal enveloping algebras}

We write $U(\lieg)$ for the universal enveloping algebra of $\lieg$ and $Z(\lieg)\subseteq U(\lieg)$ for its center. $U(\lieg)$ comes with a Lie inclusion $\lieg\to U(\lieg)$ giving rise to its universal property as the universal associative $\CC$-algebra generated by $\lieg$.

Recall the {\em Harish-Chandra map}
$$
\gamma:\;Z(\lieg)\to U(\lieh),
$$
which again by Harish-Chandra induces an isomorphism
$$
\gamma:\;Z(\lieg)\to U(\lieh)^{W(\lieh,\lieg)}
$$
which is independent of a choice of positive system.

\subsection{Infinitesimal characters}

A $(\lieg,K)$-module $X$ is said to have {\em infinitesimal character} $\lambda\in\lieh^*$, if $Z(\lieg)$ acts on $X$ via the composition
$$
Z(\lieg)\to U(\lieh)\to\CC,
$$
where the last map is given by $\lambda$. Every one-dimensional $Z(\lieg)$-module is of this form and two characters $\lambda,\lambda'\in\lieh^*$ give rise to the same infinitesimal character if and only if $w(\lambda)=\lambda'$ for some $w\in W(\lieh,\lieg)$.

\smallskip
Infinitesimal characters are very useful invariants: Every irreducible $(\lieg,K)$-module $X$ admits an infinitesimal character, in the sense that there is a unique $K$-orbit of infinitesimal characters attached to $X$. However, in our situation $K$ stabilizes $\gamma$ and therefore we need not take this into account.

\subsection{Finite-dimensional representations}

A finite-dimensional representation $V$ of $G$ is always a {\em rational representation} of the algebraic group $G$ over $\CC$. Such a $V$ is classified by Cartan's highest weight theory. Each such $V$ gives rise to a finite-dimensional {\em continuous} representation $V$ of $G(\RR)$.

\smallskip
When speaking about highest / dominant weights, we always work in the context of the Borel subalgebra $\lieb=\lieh+\lien$, i.\,e.\ dominance is understood with respect to $\lien$.

\smallskip
If $V_\lambda$ is an irreducible finite-dimensional representation of $G$ of highest weight $\lambda$, then $V_\lambda$ has infinitesimal character $\lambda+\rho$ where again $\rho=\rho(\lien)$ denotes the half sum of positive roots.

\subsection{Modules and Hecke algebras}

The category of $(\liea,B)$-modules is equivalent to the category of approximately unitial $R(\liea,B)$-modules for a suitable associative algebra $R(\liea,B)$ with approximate unit, commonly referred to as the Hecke algebra (cf.\ Chap.\ I in \cite{book_knappvogan1995}).

\smallskip
We write $\mathcal C(\liea,B)$ for the category of $(\liea,B)$-modules. It contains enough injectives and projectives.

\subsection{Induction and  production functors}

If $(\liea,B)\to(\liec,D)$ is a map of pairs, we have two forgetful functors
$$
\mathcal F:\;\mathcal C(\liec,D)\to\mathcal C(\liea,B),\quad
\mathcal F^\vee:\;\mathcal C(\liec,D)\to\mathcal C(\liea,B),
$$
with a canonical inclusion $\mathcal F\to\mathcal F^\vee$. The first one is the naive forgetful functor, while the second is given by
$$
F^\vee(X)\;=\;\Hom_{R(\liec,D)}(R(\liec,D),X)_B
$$
where the subscript $B$ denotes the subspace of $B$-finite vectors, i.\,e.\ the maximal approximately unitial submodule for the natural action of $R(\liea,B)$. Both are exact.

\smallskip
Put
\begin{align*}
P_{\liea,B}^{\liec,D}(-)\;&:=\;R(\liec,D)\otimes_{R(\liea,B)}-,\\
I_{\liea,B}^{\liec,D}(-)\;&:=\;\Hom_{R(\liea,B)}\left(R(\liec,D),-\right)_{D}.
\end{align*}
Then $P_{\liea,B}^{\liec,D}$ is left adjoint to $\mathcal F^\vee$ and $I_{\liea,B}^{\liec,D}$ is right adjoint to $\mathcal F$.

\subsection{Primary decomposition}

We say that a $(\lieg,K)$-module $X$ is {\em $Z(\lieg)$-finite} if for every $x\in X$ the space $Z(\lieg)\cdot x$ is finite-dimensional. We write $\mathcal C_{Z(\lieg){\rm -fin}}(\lieg,K)$ for the category of $Z(\lieg)$-finite modules.

\smallskip
For any character $\lambda$ of $\lieh$ the kernel of the induced algebra morphism $Z(\lieg)\to\CC$ is a maximal ideal
$$
\mathfrak m_\lambda\;:=\;\ker (Z(\lieg)\to\CC)
$$
in $Z(\lieg)$ and any maximal ideal is of this form. Therefore, every $Z(\lieg)$-finite module $X$ admits a {\em primary decomposition}
\begin{equation}
  X\;=\;\bigoplus_{\lambda\in\lieh^*/W(\lieh,\lieg)} \mathcal P_\lambda(X)
  \label{eq:primarydecomposition}
\end{equation}
where $\mathcal P_\lambda(X)$ is the the $\lambda$-primary component of $X$ (i.\,e.\ the subspace of elements annihilated by a power of $\mathfrak m_\lambda$). Each direct summand $\mathcal P_\lambda(X)$ is again a $(\lieg,K)$-module.

\smallskip
We say that $X$ has {\em generalized infinitesimal character $\lambda$}, if $X$ is $Z(\lieg)$-finite and $X=\mathcal P_\lambda(X)$. We write $\mathcal C_\lambda(\lieg,K)$ for the full subcategory of $(\lieg,K)$-modules with generalized infinitesimal character $\lambda$. In particular, we obtain a functor $\mathcal C_{Z(\lieg){\rm -fin}}(\lieg,K)\to\mathcal C_\lambda(\lieg,K)$ which by \eqref{eq:primarydecomposition} is exact.

\subsection{Orthogonality relations --- Wigner's Lemma}

Say $X$ and $Y$ are $(\lieg,K)$-modules with {\em different} generalized infinitesimal characters $\lambda$ and $\mu$, i.\,e.\ $W(\lieh,\lieg)\lambda\neq W(\lieh,\lieg)\mu$. Then
\begin{equation}
  \Ext_{\lieg,K}^\bullet(X,Y)\;=\;0.
  \label{eq:wignerslemma}
\end{equation}

\subsection{Generalities on the tensor structure}

We collect general facts about tensor products.



\subsubsection{Finite-dimensional representations}

Let $V_\lambda$, $V_\mu$ denote two irreducible finite-dimensional representations of $G$ of highest weights $\lambda$ and $\mu$ respectively. Write $\rho$ for the half sum of positive roots of $G$.

In the Grothendieck group of finite-dimensional rational representations of $G$ we then have an identity
\begin{equation}
  [V_\lambda\otimes V_\mu]\;=\;
  \sum_{\begin{subarray}c\nu\;\text{weight}\\\text{of $V_\mu$}\end{subarray}}
  {\rm m}_\mu(\nu)\cdot{\rm sgn}(\lambda+\nu+\rho)\cdot [V_{w_{\lambda+\nu+\rho}(\lambda+\nu+\rho)-\rho}]
  \label{eq:finiteexplicittensorproduct}
\end{equation}
where
\begin{itemize}
  \item ${\rm m}_\mu(\nu)$ denotes the multiplicity of $\nu$ in $V_\mu$, explicitly known by Kostant's Multiplicity Formula.
  \item $w_{\lambda+\nu+\rho}$ is any element of the Weyl group $W(\lieh,\lieg)$ such that $w_{\lambda+\nu+\rho}(\lambda+\nu+\rho)$ is dominant.
  \item ${\rm sgn}(\lambda+\nu+\rho)$ is zero if $\lambda+\nu+\rho$ is singular and the usual sign ${\rm sgn}(w_{\lambda+\nu+\rho})=(-1)^{\ell(w_{\lambda+\nu+\rho})}$ defined via the length function on $W(\lieh,\lieg)$ otherwise.
\end{itemize}
We remark that $w_{\lambda+\nu+\rho}$ is unique whenever the sign of $\lambda+\nu+\rho$ is non-zero.

\subsubsection{Tensor products with finite-dimensional representations}

Let $X$ be a $Z(\lieg)$-finite $(\lieg,K)$-module and $V$ a finite-dimensional $(\lieg,K)$-module. A theorem of Kostant (Thm.\ 7.133 in \cite{book_knappvogan1995}) guarantees that $X\otimes V$ is again $Z(\lieg)$-finite.

Furthermore, along the same lines, Kostant gives the following generalization of \eqref{eq:finiteexplicittensorproduct}: Let $X\in\mathcal C_\lambda(\lieg,K)$ and let $V_\mu$ be an irreducible rational representation of $G$ of highest weight $\lambda$. Then
\begin{equation}
  X\otimes V_\mu\;=\;\bigoplus_{\begin{subarray}c\nu\;\text{s.\,t.}\\{\rm m}_\mu(\nu)\neq 0\end{subarray}}\mathcal P_{\lambda+\nu}(X\otimes V).
  \label{eq:kostanttensorproduct}
\end{equation}

\subsubsection{Mackey Isomorphism}

If $(\liea,B)\to(\liec,D)$ is a map of pairs, we have for any $V\in\mathcal C(\liec,D)$ and $X\in\mathcal C(\liea,B)$ a natural isomorphism
\begin{equation}
  P_{\liea,B}^{\liec,D}(X\otimes\mathcal F(V))\;=\;P_{\liea,B}^{\liec,D}(X)\otimes V,
  \label{eq:mackeyisomorphism}
\end{equation}
cf.\ Theorem 2.103 in \cite{book_knappvogan1995}.

\section{Categorical setup}\label{sec:categorical}

Let $(\lieg,K)$ denote a reductive pair. We write $\mathcal C_{\rm fd}(\lieg,K)$ for the semi-simple category of finite-dimensional $(\lieg,K)$-modules and $\mathcal C_{\rm fl}(\lieg,K)$ for the category of $(\lieg,K)$-modules of finite length. Then $\mathcal C_{\rm fd}(\lieg,K)$ is symmetric monoidal with respect to the standard tensor product $\otimes$ of $\CC$-vector spaces and acts on $\mathcal C_{\rm fl}(\lieg,K)$, again via the tensor product $\otimes$ in the following sense: $\otimes$ induces the biexact functor
$$
\mathcal C_{\rm fd}(\lieg,K)\times\mathcal C_{\rm fl}(\lieg,K)\to\mathcal C_{\rm fl}(\lieg,K),\quad (X,Y)\,\mapsto\,X\otimes Y.
$$
This defines the structure of a (strong) left $\mathcal C_{\rm df}(\lieg,K)$-module on $\mathcal C_{\rm fl}(\lieg,K)$ in these sense of Yetter (cf.\ \cite{yetter1992}, see also Defintion 38 in \cite{yetter2003}).

We call a subcategory $\mathcal C$ of $\mathcal C_{\rm fl}(\lieg,K)$ {\em thick}, if for any short exact sequence
$$
\begin{CD}
  0@>>> A@>>> B@>>> C@>>> 0
\end{CD}
$$
of objects in the category $\mathcal C(\lieg,K)$ of all $(\lieg,K)$-modules the third object lies in $\mathcal C$ whenever two of the objects $A,B,C$ are objects in $\mathcal C$.

We call a thick subcategory $\mathcal C$ of $\mathcal C_{\rm fl}(\lieg,K)$ a {\em $\otimes$-submodule} over the symmetric monoidal category $\mathcal C_{\rm fd}(\lieg,K)$, if for all $X\in\mathcal C_{\rm fd}(\lieg,K)$ and all $Y\in\mathcal C$ we have $X\otimes Y\in\mathcal C$.

\smallskip
\noindent{\bf Problem:} Classify all $\otimes$-submodules of $\mathcal C_{\rm fl}(\lieg,K)$ and describe the poset of $\otimes$-submodules (wrt.\ $\subseteq$).

\smallskip
In the next section, we will execute this for $\SL_2(\RR)$.

\smallskip
We refer to \cite{stevenson2013} for a theory of supports in the context of modules for general tensor triangulated categories.

\section{The case $G=\SL(2)$}\label{sec:sl2}

For the sake of convenience, $\SU(1,1)$ is at times a good isomorphic replacement for $\SL_2(\RR)$ for the study of certain classes of representations we are interested in.

\smallskip
As $K$ we fix $\SO(2)\subseteq\SL_2(\RR)$ and obtain the reductive pair $(\liesl_2,\SO(2))$. The Cartan involution $\theta$ acts as inverse-transpose on the group level, and minus transpose on the Lie algebra.

\smallskip
Fix any Borel $B\subseteq\SL_2(\RR)$, for example
$$
B\;=\;\{\begin{pmatrix}a & x \\ 0 & a^{-1}\end{pmatrix}\},
$$
with Levi decomposition $B=TN$ where
$$
T\;=\;\{\begin{pmatrix}a & 0 \\ 0 & a^{-1}\end{pmatrix}\},\quad
N\;=\;\{\begin{pmatrix}1 & x \\ 0 & 1\end{pmatrix}\}.
$$
Then on the Lie algebra side we have
$$
\lieb\;=\;\{\begin{pmatrix}A & X \\ 0 & -A\end{pmatrix}\},
$$
and likewise, $\lieb=\liet\oplus\lien$ and
$$
\lieh\;=\;\{\begin{pmatrix}A & 0 \\ 0 & -A\end{pmatrix}\},\quad
\lien\;=\;\{\begin{pmatrix}0 & X \\ 0 & 0\end{pmatrix}\}.
$$
For compatibility with the previous section, we put $\lieh:=\liet$.

\smallskip
For the sake of concreteness, we fix the elements
\begin{align*}
H\;&:=\;\begin{pmatrix}1&0\\0&-1\end{pmatrix},\\
E\;&:=\;\begin{pmatrix}0&1\\0&0\end{pmatrix},\\
F\;&:=\;\begin{pmatrix}0&0\\1&0\end{pmatrix},
\end{align*}
which are generators of $\liet$, $\lien$ and $\lien^-$ respectively. They satisfy the relations
\begin{align*}
[E,F]\;&=\;H,\\
[H,E]\;&=\;2E,\\
[H,F]\;&=\;-2F.
\end{align*}
Recalling that the Lie bracket is altgernating, this determines the Lie algebra structure of $\liesl_2$ uniquely and may also be used to describe $\liesl_2$-modules.

\smallskip
Another conjugate choice for this triple is
\begin{align*}
H'\;&:=\;\begin{pmatrix}0&-i\\i&0\end{pmatrix},\\
E'\;&:=\;\frac{1}{2}\begin{pmatrix}1&i\\i&-1\end{pmatrix},\\
F'\;&:=\;\frac{1}{2}\begin{pmatrix}1&-i\\-i&-1\end{pmatrix}.
\end{align*}
Then $H'$ lies in the complexified Lie algebra of $K=\SO(2)\subseteq\SL_2(\RR)$, which is a convenient choice to study the action of $K$ on $(\liesl_2,\SO(2))$-modules.

\smallskip
The center $Z(\liesl_2)$ is generated by the {\em Casimir element}
$$
\Omega\;=\;H^2+1+2EF+2FE\;=\;{H'}^2+1+2E'F'+2F'E',
$$
i.\,e.\ $Z(\liesl_2)=\CC[\Omega]$ is a univariate polynomial ring.

\subsection{Finite-dimensional representations}

The standard representation $\CC^2$ of $\SL(2)$ (via matrix multiplication from the left) is an irreducible faithful representation of dimension $2$.

\smallskip
For each $m\in\NN$ there is a unique $m+1$-dimensional irreducible representation $V_\mu$ of $\SL(2)$ of $B$-highest weight
$$
\mu:\;T\to\GL(1),\quad \diag(a,a^{-1})\,\mapsto\,a^m.
$$
It is explicitly given by $V_\mu=\Sym^{m}\CC^2$. This exhausts all isomorphism classes of irreducible representations of $\SL(2)$.

\smallskip
The first standard basis vector $e_1=\begin{pmatrix}1\\0\end{pmatrix}\in\CC^2$ is annihilated by $E$ which shows that this is a highest weight vector. Likewise, $e_2=\begin{pmatrix}0\\1\end{pmatrix}\in\CC^2$ is annihilated by $F$, hence this is a lowest weight vector.

\smallskip
The action of $\liesl_2$ on $\CC^2$ is then explicitly given by matrix multiplication.

Now we may identify $\Sym^m\CC^2$ as the space of homogenous polynomials of degree $m$ in two variables, or in terms of tensor powers: The set
\begin{equation}
  \{e_1^{\otimes m-k}\otimes e_2^{\otimes k}\mid 0\leq k\leq m\}
  \label{eq:standardbasisinVmu}
\end{equation}
is a basis of $V_\mu$ consisting of generators of the weight spaces with respect to $T\subseteq\SL(2)$.

Using the explicit action of the Lie algebra on the basis vectors, the action of $Z(\liesl_2)$, i.\,e.\ of $\Omega$ on $V_\mu$ is easily computed:
\begin{align*}
  \Omega\cdot e_1^{\otimes m}\;
  &=\;H^2\cdot e_1^{\otimes m}+1\cdot e_1^{\otimes m}+2EF\cdot e_1^{\otimes m}+\underbrace{2FE\cdot e_1^{\otimes m}}_{=\,0}\\
  &=\;m^2\cdot e_1^{\otimes m}+e_1^{\otimes m}+2m\cdot E\cdot e_1^{\otimes m-1}\otimes e_2\\
  &=\;(m+1)^2\cdot e_1^{\otimes m}.
\end{align*}
This result is not surprising, because under the Harish-Chandra isomorphism
\begin{equation}
  \gamma(\Omega)\;=\;H^2.
  \label{eq:harishchandracasimir}
\end{equation}

\smallskip
We remark that in this case the Casimir operator distinguishes isomorphism classes of finite-dimensional representations: The infinitesimal character of $V_\mu$ is $\mu+\rho$ in Harish-Chandra's parametrization.

\subsection{Principal series representations}

In this subsection we recall well-known statements about principal series representations and subsequently study tensor products of principal series representations with finite-dimensional representations. The technical main results are Propositions \ref{prop:irreduciblecase} and \ref{prop:reduciblecase}.

\subsubsection{Parabolic induction and principal series}

Recall the {\em Iwasawa decomposition}
\begin{equation}
  \SL_2(\RR)\;=\;B\cdot K\;=\;N\cdot T\cdot\SO(2)
  \label{eq:iwasawa}
\end{equation}
which yields the following unique representation of any $g\in\SL_2(\RR)$:
$$
g\;=\;
\begin{pmatrix}
  \sqrt{y} & x/\sqrt{y} \\ 0 & 1/\sqrt{y}
\end{pmatrix}\cdot
\begin{pmatrix}
  \cos\alpha & -\sin\alpha \\ \sin\alpha & \cos\alpha
\end{pmatrix}
$$
with $x\in\RR$, $y>0$, $0\leq\alpha<2\pi$. We remark that the image of $i$ under the action of $g$ and hence also under the action of the first matrix right hand side on the upper half plane is $x+yi$.

\smallskip
In this representation of $g$, we only take the connected component $A$ of the identity in $T(\RR)$ into account, because
$$
M\;=\;T(\RR)\cap\SO(2)\;=\;\{\pm{\bf1}_2\},
$$
which yields the {\em Langlands decomposition}:
$$
\SL_2(\RR)\;=\;N\cdot A\cdot M\cdot K.
$$
Therefore, describing a character of $T(\RR)$ is the same as a character of $A$ together with a character of $M$.

\smallskip
Characters of $A\cong\RR$ are parametrized by complex numbers $\lambda\in\CC$, to each of which we attach the character
$$
|\,\cdot\,|^{\lambda}:\;A\to\CC^\times,\quad\diag(a,a^{-1})\,\mapsto\,a^{\lambda}.
$$
In the special case of $\rho$ half sum of roots in $N$ we get $|\,\cdot\,|^{\rho}=|\,\cdot\,|^1$ since $\rho=1$. This is the square root of the modular character $\Delta_B$ restricted to $A$.

\smallskip
Characters of $M$ are indexed by parity $\epsilon\in\{0,1\}$, $0$ representing the trivial character and $1$ representing the {\em sign character} ${\rm sgn}:M\to\{\pm1\}$. In summary, $\epsilon$ corresponds to the character ${\rm sgn}^\epsilon$.

\smallskip
Hence, we parametrized characters $\chi$ of $T(\RR)=AM$ by pairs $(\lambda,\epsilon)\in\CC\times\{0,1\}$.

\smallskip
Pulling back these characters to $B$ we obtain the {\em principal series representation} attached to $(\lambda,\epsilon)$ as
\begin{align*}
  I_{\lambda,\epsilon}\;&:=\;\Ind_{B}^{\SL_2(\RR)}{\rm sgn}^\epsilon|\,\cdot\,|^{\lambda+\rho}\\
  &:=\;\{f\in\mathcal C^\infty(\SL_2(\RR))\mid \forall b\in B,g\in\SL_2(\RR):f(b\cdot g)={\rm sgn}^\epsilon(b)|b|^{\lambda+\rho}f(g)\,\wedge\,f\;\text{right $K$-finite}\}.
\end{align*}
This is a $(\liesl_2,\SO(2))$-module and as such
\begin{equation}
  I_{\lambda,\epsilon}\;=\;P_{\lieb,M}^{\liesl_2,\SO(2)}\left({\rm sgn}^\epsilon|\,\cdot\,|^{\lambda-\rho}\right)\;=\;R(\liesl_2,\SO(2))\otimes_{U(\liesl_2),M}U(\liesl_2)\otimes_{U(\lieb)}{\rm sgn}^\epsilon|\,\cdot\,|^{\lambda-\rho}.
  \label{eq:bernsteinfunctordescription}
\end{equation}
  Its infinitesimal character is $\lambda$: We let $\Omega$ act on $1\otimes 1$ on the right hand side and obtain\footnote{The application of the Bernstein functor does not change infinitesimal characters.}
  \begin{align*}
    \Omega(1\otimes 1)\;
    &=\;H^2\otimes1+2H\otimes 1+1\otimes1+\underbrace{4FE\otimes 1}_{=\,0}\\
    &=\;\left((\lambda-\rho)(H)^2+2(\lambda-\rho)(H)+1\right)\cdot(1\otimes 1)\\
    &=\;\lambda^2\cdot(1\otimes1)
  \end{align*}
  which is consistent with our claim thanks to \eqref{eq:harishchandracasimir}.

\subsubsection{The compact picture}

The Iwasawa decomposition \eqref{eq:iwasawa} may be exploited to study the $K$-action on $I_{\lambda,\epsilon}$. As a $K$-module, we have
  \begin{align*}
    I_{\lambda,\epsilon}|_K\;&=\;\Ind_{B\cap K}^{K}{\rm sgn}^\epsilon|\,\cdot\,|^{\lambda+\rho}\\
    &=\;\Ind_{M}^{K}{\rm sgn}^\epsilon|\,\cdot\,|^{\lambda+\rho}\\
    &=\;\Ind_{M}^{K}{\rm sgn}^\epsilon\\
    &=\;\{f\in\mathcal C^\infty(\SO(2))\mid f(- g)=(-1)^\epsilon f(g)\,\wedge\,f\;\text{right $K$-finite}\}.
  \end{align*}
  By classical Fourier theory, the $K$-finite smooth functions on the unit circle are given by the finite linear combinations of the characters
  $$
  \chi_k:\;\SO(2)\to\CC^\times,\quad \begin{pmatrix}
  \cos\alpha & -\sin\alpha \\ \sin\alpha & \cos\alpha
\end{pmatrix}\,\mapsto\,e^{i k \alpha}
  $$
  for $k\in\ZZ$. Hence, the above description shows:
  \begin{equation}
  I_{\lambda,\epsilon}|_K\;=\;\bigoplus_{k\equiv\epsilon\,(2)}\chi_k.
  \end{equation}

  \smallskip
  Therefore, similar to the basis \eqref{eq:standardbasisinVmu} in the finite-dimensional case, $I_{\lambda,\epsilon}$ admits a basis $\{w_k\}_{k\equiv \epsilon\,(2)}$ consisting of $\SO(2)$-eigenfunctions, numbered by integers with same parity as $\epsilon$. This basis may be normalized in such a way that
  \begin{align*}
    H'\cdot w_k\;&=\;k\cdot w_k,\\
    E'\cdot w_k\;&=\;\frac{\lambda+k+1}{2}\cdot w_{k+2},\\
    F'\cdot w_k\;&=\;\frac{\lambda-k+1}{2}\cdot w_{k-2}.
  \end{align*}
  In fact, comparing those relations with the finite-dimensional explains all cases where $I_{\lambda,\epsilon}$ is reducible with one exception.

\subsubsection{Decomposition of principal series}

The following facts about $I_{\lambda,\epsilon}$ are well known:
\begin{itemize}
  \item The dual of a principal series is a principal series: $I_{\lambda,\epsilon}^\vee\cong I_{-\lambda,\epsilon}$.
  \item If $I_{\lambda,\epsilon}$ is irreducible, it is self-dual, i.\,e.\ $I_{\lambda,\epsilon}\cong I_{-\lambda,\epsilon}$.
  \item $I_{\lambda,\epsilon}$ is reducible if and only if $\lambda\in\ZZ$ and $\lambda\not\equiv\epsilon\pmod{2}$.
  \item $I_{0,1}=D_{0}^+\oplus D_{0}^-$ decomposes into the holomorphic and anti-holomorphic limits of discrete series representations, where $D_0^{\pm}$ is generated by $w_{\pm k}$ for $k\geq 1$.
  \item If $0<\lambda\in\ZZ$ and $\lambda\not\equiv\epsilon\pmod{2}$, then we have a non-split short exact sequence
    \begin{equation}
      0\to D_{\lambda}^+\oplus D_{\lambda}^-\to I_{\lambda,\epsilon}\to V_{\lambda-1}\to 0
      \label{eq:positivedecomposition}
    \end{equation}
    where $D_{\lambda}^\pm$ denotes the $\lambda$-th holomorphic and anti-holomorphic discrete series representations, generated by $w_{\pm k}$ for $k\geq\lambda+1$, i.\,e.\ $\chi_{\pm(\lambda+1)}$ is the minimal $K$-type in $D_{\lambda}^\pm$.
  \item If $0>\lambda\in\ZZ$ and $\lambda\not\equiv\epsilon\pmod{2}$, then we have a non-split short exact sequence
    \begin{equation}
    0\to V_{\lambda-1}
    \to I_{\lambda,\epsilon}
    \to D_{\lambda}^+\oplus D_{\lambda}^-
    \to 0.
      \label{eq:negativedecomposition}
    \end{equation}
\end{itemize}

The sequence \eqref{eq:positivedecomposition} and its dual \eqref{eq:negativedecomposition} are manifestations of the fact that
$$
H^1(\liesl_2,\SO(2); D_{\lambda}^\pm\otimes V_{\lambda-1})\;=\;\Ext^1_{\liesl_2,\SO(2)}(D_{\lambda}^\mp,V_{\lambda-1})\;=\;\CC,
$$
which explains the relevance of the discrete series to modular forms from a topological point of view.

\subsubsection{Classification of irreducible $(\liesl_2,\SO(2))$-modules}

By the Casselman-Wallach embedding theorem, each irreducible admissible representation of $G$ embeds into a principal series representation. Therefore, we deduce from the above description of all principal series representations for $\SL_2(\RR)$, a classification of irreducible $(\liesl_2,\SO(2))$-modules\footnote{Strictly speaking we also invoke a globalization result here: Each Harish-Chandra module admits a globalization, i.\,e.\ extends to a representation of $G$.}.

\subsubsection{Tensor products with finite-dimensional representations}

The Mackey isomorphism \eqref{eq:mackeyisomorphism} together with the intrinsic construction \eqref{eq:bernsteinfunctordescription} of $I_{\lambda,\epsilon}$ implies that for any finite-dimensional irreducible representation $V_\mu$ of highest weight $\mu$, we have
$$
I_{\lambda,\epsilon}\otimes V_\mu\;=\;
P_{\lieb,M}^{\liesl_2,\SO(2)}\left({\rm sgn}^\epsilon|\,\cdot\,|^{\lambda-\rho}\otimes V_\mu\right)\;=\;
\Ind_{B\cap K}^{K}\left({\rm sgn}^\epsilon|\,\cdot\,|^{\lambda+\rho}\otimes V_\mu\right).
$$
Now as $B$-module, $V_\mu$ admits a filtration
$$
0\subsetneq
V_\mu^{0}\subsetneq
V_\mu^{1}\subsetneq
V_\mu^{2}\subsetneq
\dots\subsetneq
V_\mu^{\mu-2}\subsetneq
V_\mu^{\mu-1}\subsetneq
V_\mu^{\mu}=V_\mu
$$
whose successive quotients are the one-dimensional weight spaces:
$$
V_\mu^{j}/
V_\mu^{j-1}\;\cong\;\CC\cdot e_1^{\otimes(\mu-j)}\otimes e_2^{\otimes j}.
$$
By the exactness and faithfulness of the induction functor, we obtain a filtration
\begin{align*}
0\;&\subsetneq\;
\Ind_{B\cap K}^{K}{\rm sgn}^\epsilon|\,\cdot\,|^{\lambda+\rho}\otimes V_\mu^{0}\\
&\subsetneq\;
\Ind_{B\cap K}^{K}{\rm sgn}^\epsilon|\,\cdot\,|^{\lambda+\rho}\otimes V_\mu^{1}\\
&\subsetneq\;
\dots\\
&\subsetneq\;
\Ind_{B\cap K}^{K}{\rm sgn}^\epsilon|\,\cdot\,|^{\lambda+\rho}\otimes V_\mu^{\mu-1}\\
&\subsetneq\;
\Ind_{B\cap K}^{K}{\rm sgn}^\epsilon|\,\cdot\,|^{\lambda+\rho}\otimes V_\mu^{\mu}=I_{\lambda,\epsilon}\otimes V_\mu
\end{align*}
of $I_{\lambda,\epsilon}\otimes V_\mu$, with successive quotients
\begin{equation}
  \Ind_{B\cap K}^{K}\left({\rm sgn}^\epsilon|\,\cdot\,|^{\lambda+\rho}\otimes V_\mu^{j}\right)/
  \Ind_{B\cap K}^{K}\left({\rm sgn}^\epsilon|\,\cdot\,|^{\lambda+\rho}\otimes V_\mu^{j-1}\right)\;\cong\;
  I_{\lambda+\mu-2j,\epsilon+\mu\,(2)}.
\end{equation}
Note that the non-trivial element $-{\bf1}_2$ of the center acts on $V_\mu^{j}/V_\mu^{j-1}$ as $(-1)^\mu$.

In order to refine this filtration, we distinguish two cases:
\begin{itemize}
\item[(a)] $I_{\lambda,\epsilon}$ is irreducible. Then either
  \begin{itemize}
    \item[(${\rm a}_1$)] $\lambda\not\in\ZZ$ in which case the same is true for all $\lambda+\mu-2j$, which implies that all composition factors $I_{\lambda+\mu-2j,\epsilon+\mu\,(2)}$ are again irreducible.
    \item[(${\rm a}_2$)] $\lambda\in\ZZ$ and $\lambda\equiv\epsilon\pmod{2}$. Then
      $$
      \lambda+\mu-2j\equiv \epsilon+\mu\,\pmod{2}
      $$
      again, hence the composition factors $I_{\lambda+\mu-2j,\epsilon+\mu\,(2)}$ are once again irreducible.
  \end{itemize}
\item[(b)] $I_{\lambda,\epsilon}$ is reducible. Then, by the same computation as in case ${\rm a}_2$, $I_{\lambda+\mu-2j,\epsilon+\mu\,(2)}$ is reducible as well.
\end{itemize}

In case (a), we have
$$
I_{\lambda+\mu-2j,\epsilon+\mu\,(2)}\;=\;I_{-\lambda-\mu+2j|,\epsilon+\mu\,(2)}
$$
by self-duality, which leads to the question if the relation
$$
\lambda+\mu-2j\;=\;-\lambda-\mu+2j'
$$
can be satisfied. In fact, this amounts to
$$
\lambda+\mu\;=\;j+j',
$$
i.\,e.\ this is only possible if $\lambda\in\ZZ$ and since $0\leq j,j'\leq \mu$, we get different $j,j'$ satisfying this relation if and only if
\begin{equation}
\mu-1\,\geq\,\lambda\,\geq\,1-\mu
\label{eq:lambdamubound}
\end{equation}
We conclude that in case (a), if $\lambda\not\in\ZZ$ or if \eqref{eq:lambdamubound} is violated, the infinitesimal characters of all composition factors are pairwise distinct, which by Wigner's Lemma then implies
\begin{proposition}[Irreducible principle series]\label{prop:irreduciblecase}
  In case (a), if $\lambda\not\in\ZZ$ or if \eqref{eq:lambdamubound} is not satisfied,
  \begin{equation}
    I_{\lambda,\epsilon}\otimes V_\mu\;=\;\bigoplus_{j=0}^\mu I_{\lambda+\mu-2j,\epsilon+\mu\,(2)}.
  \end{equation}
  In case (a), if $\lambda\in\ZZ$ and condition \eqref{eq:lambdamubound} is satisfied,
  \begin{equation}
    I_{\lambda,\epsilon}\otimes V_\mu\;=\;
    \bigoplus_{j=0}^{\lambda-1}
    I_{\lambda+\mu-2j,\epsilon+\mu\,(2)}
    \,
    \oplus
    \,
    \bigoplus_{j=\lambda}^{\min(\mu-1,\lambda+\mu-1)}
    \mathcal{I}_{\lambda-\mu+2j,\epsilon+\mu\,(2)}
    \,
    \oplus
    \,
    \bigoplus_{j=\lambda+\mu+1}^{\mu}
    I_{\lambda+\mu-2j,\epsilon+\mu\,(2)}
  \end{equation}
  where $\mathcal{I}_{\lambda-\mu+2j,\epsilon+\mu\,(2)}$ is of length two and fits into a short exact sequence
  \begin{equation}
    0\to
    I_{\lambda+\mu-2j,\epsilon+\mu\,(2)}\to
    \mathcal{I}_{\lambda-\mu+2j,\epsilon+\mu\,(2)}\to
    I_{\lambda+\mu-2j,\epsilon+\mu\,(2)}\to
    0.
  \end{equation}
\end{proposition}

More generally, in case (b), condition \eqref{eq:lambdamubound} characterizes the situations, where two composition factors share the same infinitesimal character. Therefore, we obtain along the same lines
\begin{proposition}[Reducible principle series]\label{prop:reduciblecase}
  In case (b), if \eqref{eq:lambdamubound} is not satisfied,
  \begin{equation}
    I_{\lambda,\epsilon}\otimes V_\mu\;=\;\bigoplus_{j=0}^\mu I_{\lambda+\mu-2j,\epsilon+\mu\,(2)}.
  \end{equation}
  In case (b), if condition \eqref{eq:lambdamubound} is satisfied,
  \begin{equation}
    I_{\lambda,\epsilon}\otimes V_\mu\;=\;
    \bigoplus_{j=0}^{\lambda-1}
    I_{\lambda+\mu-2j,\epsilon+\mu\,(2)}
    \,
    \oplus
    \,
    \bigoplus_{j=\lambda}^{\min(\mu-1,\lambda+\mu-1)}
    \mathcal{I}_{\lambda-\mu+2j,\epsilon+\mu\,(2)}
    \,
    \oplus
    \,
    \bigoplus_{j=\lambda+\mu+1}^{\mu}
    I_{\lambda+\mu-2j,\epsilon+\mu\,(2)}
  \end{equation}
  where $\mathcal{I}_{\lambda-\mu+2j,\epsilon+\mu\,(2)}$ is of length two and fits into a short exact sequence
  \begin{equation}
    0\to
    I_{\lambda+\mu-2j,\epsilon+\mu\,(2)}\to
    \mathcal{I}_{\lambda-\mu+2j,\epsilon+\mu\,(2)}\to
    I_{-\lambda-\mu+2j,\epsilon+\mu\,(2)}\to
    0.
  \end{equation}
\end{proposition}

\subsection{Asymptotic $K$-support}

Before we can fully decompose tensor products of irreducible with finite-dimensional representations, we need the notion of {\em asymptotic $K$-support}.

\smallskip
Let $S\subseteq\RR^n$ denote a subset. The {\em asymptotic cone} $S\infty$ attached to $S$ is the closed cone
$$
S\infty\;:=\;\{y\in\RR^n\mid\exists\;\text{a sequence}\;(y_n,\varepsilon_n)_{n\in\NN}\in S\times\RR_{>0}:\;\lim \varepsilon_ny_n=y\,\wedge\,\lim\varepsilon_n=0\}.
$$

The asmptotic $K$-support $\AS_K(X)$ of a $(\lieg,K)$-module $X$ is defined as the asymptotic cone attached to the lattice of highest weights of the $K$-types occuring in $X$.

\smallskip
The asymptotic $K$-support is of interest of us, because for any finite-dimensional representation $V$ we have $$
\AS_K(X\otimes V)\;=\;\AS_K(X).
$$

In our case at hand, the situation is the following:
\begin{table}[h!]
  \begin{center}
    \label{tab:table1}
    \begin{tabular}{l|c|c} 
      $X$ & \textbf{$K$-types} & $\AS_K(X)$\\
      \hline
      $V_\mu$ & $\{-\mu,\mu+2,\dots,\mu-2,\mu\}$ & $\{0\}$\\
      $I_{\lambda,\epsilon}$ & $\epsilon+2\ZZ$ & $\RR$\\
      $D_\ell^+$ & $\{\ell+1,\ell+3,\dots\}$ & $\RR_{\geq 0}$\\
      $D_\ell^-$ & $\{-\ell-1,-\ell-3,\dots\}$ & $\RR_{\leq 0}$
    \end{tabular}
    \caption{Asymptotic $K$-support for irreducible $(\liesl_2,\SO(2))$-modules.}
  \end{center}
\end{table}

\subsection{Applications to tensor products with irreducibles}

In Proposition \ref{prop:irreduciblecase} we treated tensor products of irreducible principle series with finite-dimensional representations. In this section, we treat the remaining cases: tensor products of finite-dimensionals and tensor products of (limits of) discrete series representations with finite-dimensional representations.

\smallskip
In our situation, \eqref{eq:finiteexplicittensorproduct} specializes to
\begin{proposition}[Finite-dimensional]\label{prop:finitecase}
  Let $\lambda,\mu\geq0$. Then
  \begin{equation}
    V_\lambda\otimes V_\mu\;=\;\bigoplus_{j=0}^{\min(\lambda,\mu)} V_{\lambda+\mu-2j}
    \label{eq:finitetensorproducts}
  \end{equation}
\end{proposition}

\smallskip
The invariance of the asymptotic $K$-support tells us that $D_\ell^\pm\otimes V_\mu$ only admits finite-dimensional representations and (limits of) discrete series representations with the same $K$-support as subquotients.

\smallskip
Now $D_{\ell}^\pm$ is a subquotient of $I_{\ell,\ell+1\,(2)}$, hence by the exactness of the tensor product $D_{\ell}^\pm\otimes V_\mu$ is a subquotient of $I_{\ell,\ell+1\,(2)}\otimes V_\mu$.

\smallskip
By Proposition \ref{prop:reduciblecase} we know which composition factors can occur in $I_{\ell,\ell+1\,(2)}\otimes V_\mu$, and since for all composition factors occuring in $I_{\ell,\ell+1\,(2)}\otimes V_\mu$ need to be accounted for by one of the three cases representations $V_{\ell-1}\otimes V_\mu$ (finite-dimensional), $D_\ell^+\otimes V_\mu$ and $D_\ell^-\otimes V_\mu$, we deduce:


\begin{proposition}[Discrete series]\label{prop:discretecase}
  Let $\ell\geq 0$ and $\mu\geq 0$. We consider the condition

  Then, if $\ell\geq\mu$,
  \begin{equation}
    D_\ell^\pm\otimes V_\mu\;=\;
    \bigoplus_{j=0}^{\mu} D_{\ell+\mu-2j}^\pm.
  \end{equation}
  Otherwise, if $\ell < \mu$,
  \begin{equation}
    D_\ell^\pm\otimes V_\mu\;=\;
    \bigoplus_{j=0}^{\ell-1} D_{\ell+\mu-2j}^\pm
    \,\oplus\,
    \bigoplus_{j=\ell}^\mu \mathcal D_{\ell+\mu-2j}^\pm
  \end{equation}
  where $\mathcal{D}_{\ell-\mu+2j,\epsilon+\mu\,(2)}^\pm$ is of length three (two if $\mu=\ell+2j$) and admits $V_{\mu-\ell-2j-1}$ and $\mathcal D_{\ell+\mu-2j}^\pm$ as composition factors, where the (limits of) discrete series occurs twice.
\end{proposition}

\subsection{Classification of thick $\otimes$-submodules}

Summing up, we obtain the following classification of thick $\mathcal C_{\rm fd}(\SL(2))$-submodules of the category $\mathcal C_{\rm fl}(\liesl_2,\SO(2))$ of finite length $(\liesl_2,\SO(2))$-modules.

\begin{theorem}\label{thm:classification}
  Any thick $\mathcal C_{\rm fd}(\SL(2))$-submodule of the category $\mathcal C_{\rm fl}(\liesl_2,\SO(2))$ of finite length $(\liesl_2,\SO(2))$-modules is of the following form:
\begin{itemize}
  \item[(i)] The zero subcategory.
  \item[(ii)] $\mathcal C_{\rm fd}(\SL(2))$.
  \item[(iii)] $\mathcal C^\pm$: the subcategory of modules whose composition factors are either finite-dimensional or a (limits of) holomorphic (resp.\ anti-holomorphic) discrete series.
  \item[(iv)] $\mathcal C(\lambda,\epsilon)$ for $\lambda\in\CC$ and $\epsilon\in\{0,1\}$ with either $\lambda\not\in\ZZ$ or $\lambda\not\equiv\epsilon\pmod{2}$: The subcategory of modules whose composition factors are irreducible principle series $I_{\lambda+j,\epsilon+j\,(2)}$ for any $j\in\ZZ$.
  \item[(v)] Any combination of the above.
\end{itemize}
\end{theorem}

\begin{remark}
  In (iv) we have $\mathcal C(\lambda,\epsilon)=\mathcal C(\lambda',\epsilon')$ if and  only if $\lambda-\lambda'\in\ZZ$ and $\epsilon+\lambda\equiv\epsilon'+\lambda'\pmod{2}$.
\end{remark}

\begin{remark}
  We remark that the categories in (ii) and (iv) are minimal $\otimes$-submodules. The category $\mathcal C^\pm$ in (iii) always contains the category $\mathcal C_{\rm fd}(\SL(2))$ of finite-dimensional modules.
\end{remark}

\begin{remark}\label{rmk:Xirred}
  Each $\otimes$-submodule $\langle X\rangle$ generated by an irreducible $(\liesl_2,\SO(2))$-module $X$ belongs to one of the three cases (ii), (iii) and (iv) and each $\otimes$-submodule listed in (ii), (iii) and (iv) is generated by an irreducible $(\liesl_2,\SO(2))$-module, the latter not being unique.
\end{remark}

\begin{remark}\label{rmk:Xall}
  The statement in (v) signifies that for every subset $S$ of cases listed in (ii), (iii) and (iv) we have a unique $\otimes$-submodule $\mathcal C(S)$ consisting of those Harish-Chandra modules, whose composition factors belong to the cases belonging to $S$. If we restrict to the set $\mathfrak{S}$ of subsets $S$ with the property, that if at least one case from (iii) belongs to $S$, then also (ii) belongs to $S$, then $S\mapsto\mathcal C(S)$ induces a bijection
  $$
  \mathfrak{S}\,\to\,\{\otimes\text{-submodules of }\mathcal C_{\rm fl}(\liesl_2,\SO(2))\}
  $$
  which is an isomorphism of posets (both sides wrt.\ to \lq{}$\subseteq$\rq).

  This picture may be topologized as follows. Consider the set $\mathfrak{C}$ of all non-trivial cases, i.\,e.\ the set of all cases occuring in (ii), (iii) and (iv) and topologize this set as follows: Each singleton corresponding to case (ii) or one case in (iv) is open. Furthermore, if $x_{(ii)}$ represents case (ii), then for each case $y$ in (iii), the set $\{x_{(ii)},y\}$ is open. Then $\mathfrak{S}$ is the set of closed subsets in $\mathfrak{C}$.
\end{remark}

\begin{proof}[{Proof of Theorem \ref{thm:classification}}]
  Let $\mathcal C$ denote a thick $\otimes$-submodule of $\mathcal C_{\rm fl}(\liesl_2,\SO(2))$. If $\mathcal C\neq 0$, it contains an irreducible module $X$. If $X$ is finite-dimensional we have $\mathcal C_{\rm fd}(\SL(2))\subseteq\mathcal C$. By Proposition \ref{prop:discretecase} the same applies if $X$ is a discrete series representation or a limits of discrete series representation. Furthermore, by the same proposition, depending on whether $X$ belongs to the holomorphic or an anti-holomorphic discrete series (or is a limit thereof), $\mathcal C^+$ or $\mathcal C^-$ are contained in $\mathcal C$.

  \smallskip
  Likewise, if $X=I_{\lambda,\epsilon}$ is an irreducible principal series representation, Proposition \ref{prop:irreduciblecase} implies that $\mathcal C(\lambda,\epsilon)\subseteq\mathcal C$.

  \smallskip
  Repeating this argument ad infinitum (transfinite induction), we exhaust all composition factors occuring in objects of $\mathcal C$ with the cases (ii), (iii) and (iv) in the theorem.
\end{proof}

\subsection{Geometric interpretation}\label{sec:geometricinterpretation}

There is a close analogy between the classification in Theorem \ref{thm:classification} and $\SO(2)$-orbits on the flag variety $\PP^1=G/B$ over $\CC$ in the Beilinson-Bernstein classification:
\begin{itemize}
  \item[(ii)] The category $\mathcal C_{\rm fd}(\SL(2))$ corresponds to the closed orbit $\PP^1(\CC)$ for the compact form $\SU(2)$ of $\SL_2(\RR)$.
  \item[(iii)] The categories $\mathcal C^\pm$ correspond to the two closed $\SO(2)$-orbits given by the poles $\{0\}$ and $\{\infty\}$ in $\PP^1(\CC)$, which in turn are contained in the closed orbit for the compact form $\SU(2)$ of $\SL_2(\RR)$.
   \item[(iv)] $\mathcal C(\lambda+\ZZ,\epsilon)$ for $\lambda\in\CC$ and $\epsilon\in\{0,1\}$ with either $\lambda\not\in\ZZ$ or $\lambda\not\equiv\epsilon\pmod{2}$ correspond to the open orbit $\CC^\times\subseteq\PP^1(\CC)$.
\end{itemize}

We remark that while the irreducible objects in $\mathcal C_{\rm fd}(\SL(2))$ except for the tensor unit are not unitary for $\SL_2(\RR)$, they are unitary representations of $\SU(2)$. Bernstein offers in \cite{book_bernstein2018} an approach to resolve this.

\subsection{Classification of irreducible modules}\label{sec:completeclassification}

The $\otimes$-submodules in (ii), (iii) and (iv) all contain representations with various infinitesimal characters. Consequently, we may think of $\Spec Z(\liesl_2)$ geometrically as a stratum transversal to the classification in Theorem \ref{thm:classification}. Putting these two data together, i.\,e.\ by identifying the infinitesimal characters appearing in in each of the $\otimes$-submodules in (ii), (iii) and (iv), we may deduce a classification of irreducible $(\lieg,K)$-modules as follows.

Following remark \ref{rmk:Xirred}, we have a map
$$
\chi:\;\{\,{\rm irreducible }(\liesl_2,\SO(2))\text{-modules}\,\}/\!\!\cong\;\to\,{\mathfrak{S}},\quad X\,\mapsto\,\overline{\left\{\langle X\rangle\right\}},
$$
sending an isomorphim class of an irreducible $(\liesl_2,\SO(2))$-module $X$ to the closure of the singleton containing the $\otimes$-submodule generated by it. We remark that by the definition of the topology on $\mathfrak{C}$, the image of $\chi$ consists of the (topologically) {\em irreducible} subsets of $\mathfrak{C}$. Different elements $X,Y$ in a fiber $\chi^{-1}(x)$ have different infinitesimal characters $\lambda(X)$, $\lambda(Y)$. Hence the pair $(\chi(X),\lambda(X))$ determines the isomorphism class of $X$ uniquely. However, the infinitesimal characters of different elements in the same fiber differ by an {\em integer}. Hence, fixing base-point $\lambda_{(\lambda+\ZZ,\epsilon)}\in\CC$ for the set of infinitesimal characters occuring in each $\mathcal C(\lambda+\ZZ,\epsilon)$, we obtain a bijection
$$
\{\,{\rm irreducible }(\liesl_2,\SO(2))\text{-mod.}\,\}/\!\!\cong\;\to\,\{(C,\lambda)\in{\mathfrak{S}}\times\ZZ\mid C\;\text{irred.\ and}\;\mathcal C_{\rm fd}(\SL(2))\in C\Rightarrow\lambda\geq 1\},
$$
explicitly given by
$$
X\;\mapsto\;
\begin{cases}
  \left(\chi(X),\lambda(X)\right),&\text{$X$ finite-dimensional or a (limits of) discrete series,}\\
  \left(\chi(X),\lambda(X)-\lambda_{(\lambda(X)+\ZZ,\epsilon)}\right),&\text{$X$ an irreducible principal series.}
\end{cases}
$$

\addcontentsline{toc}{section}{References}
\bibliographystyle{plain}


\end{document}